\documentclass[final,onefignum,onetabnum]{siamart171218}

\usepackage{tikz}
\usepackage{amsfonts}
\usepackage{amssymb}
\usepackage{graphicx}
\usepackage{epstopdf,amsmath,stmaryrd}
\usepackage{algorithmic}
\usepackage{amsopn}

\ifpdf
  \DeclareGraphicsExtensions{.eps,.pdf,.png,.jpg}
\else
  \DeclareGraphicsExtensions{.eps}
\fi

\newsiamremark{remark}{Remark}
\newsiamremark{hypothesis}{Hypothesis}
\crefname{hypothesis}{Hypothesis}{Hypotheses}
\newsiamthm{claim}{Claim}

\newcommand{\pd}[3][1]{
  \ifnum#1=1
  \frac{\partial{#2}}{\partial{#3}}
  \else
  \frac{\partial^{#1}{#2}}{\partial{{#3}^{#1}}}
  \fi
 }

 \newcommand{\sgn}{\mbox{sgn}}

 \newcommand{\od}[3][1]{
  \ifnum#1=1
  \frac{d{#2}}{d{#3}}
  \else
  \frac{d^{#1}{#2}}{d{#3}^{#1}}
  \fi
}

\begin{document}

\headers{Graphene Plasmons}{F.~Santosa and T.~Shi}

\title{\bf Analysis and Simulation of Plasmons in Graphene with Time- and Space-dependent Drude Weight}\thanks{Submitted to the editors DATE.}

\author{Fadil Santosa\thanks{Applied Mathematics and Statistics, Johns Hopkins University (\email{fsantos9@jhu.edu})}
\and Tong Shi\thanks{School of Mathematics, University of Minnesota
  (\email{shixx739@umn.edu}).}
}





\maketitle

\begin{abstract}
We study the propagation of plasmons on graphene. The problem is considered in two dimensions with a transverse magnetic (TM) electromagnetic field. The graphene material is assumed to be flat and is modeled as a conductive sheet. This leads to a jump condition for the magnetic field on the sheet where it is related to the current density on the sheet. The current density itself satisfies Drude's law. The model then consists of Maxwell's equations coupled to current density on the sheet. To make the problem more computationally approachable, we develop a time-dependent integro-differential equation for the current density. This effectively reduces the problem to one space dimension. A finite difference method is proposed to solve the resulting equation. Numerical examples are given to illustrate previously unreported behavior of the system. 
\end{abstract}



\section{Introduction} 
\label{Section 1}

Two-dimensional materials such as graphene, which is formed by lattice of carbon atoms, have attracted great interest in the study of plasmons Jablan et al,  \cite{Plasmon_Jablan}, Low et al \cite{Low_2017}.
 Compared to traditional metal-dielectric materials, this highly conductive material has certain advantages. It is known that the carrier concentration of the graphene sheet can be tunable Chen et al \cite{Graphene_Tunable}.
 Additionally, graphene is capable of producing of low-loss plasmons Garcia de Abajo \cite{Graphene_Challenge}. 
Comprehensive review of this area can be found in Goncalves and Peres \cite{Intro_Book} and Fan \cite{Intro_Fan}.

There are many applications of graphene plasmons. For example, optical sensors based on graphene and other two-dimensional materials have been widely studied; see Zhang et al \cite{Sensor_Zhang}.
Photodetectors based on graphene, which is a sensor of electromagnetic radiation and converts the photon energy into electrical signal, has been investigated and developed in Koppens et al \cite{Photodetector_Koppens},
and Wang et al \cite{Photodetector_Wang}. Also, graphene plasmons can be used to confine and manipulate terahertz waves, which enables improved performance and functionality for terahertz devices as described in Low and Avouris \cite{Low_2014}, and Castilla et al \cite{Terahertz_Castilla}. More recently, graphene-based tactile sensors \cite{kim-tactile}, gas sensors \cite{khaliji-gas}, and spectrometer \cite{lee-spectroscopy} have been proposed and developed.

The present work is a follow-up on a work started in Wilson et al \cite{Josh_Physics_Review} and \cite{Maxwell_Fadil}. In those papers, the authors studied plasmon propagation in graphene whose Drude weight is time-dependent. Among the phenomena observed and analyzed are reversal and amplification of plasmons through modulation of the Drude weight as a function of time. In \cite{Maxwell_Fadil} the authors derived a single integro-differential equation which governs the amplitude of the current density on the graphene sheet. This allows for a rigorous analysis of the plasmon model and the development of an accurate numerical scheme for simulation.

In this work, we continue the development began in \cite{Maxwell_Fadil} and study the current density evolution on a graphene sheet with time- and space-dependent Drude weight. This requires a generalization of the previous results. Instead of a single integro-differential equation for amplitude, we now must consider a general form for the current density. This leads to a new equation for current density which has the form of a partial integro-differential equation. Little is known about this new equation as we have not seen anything similar in the open literature. Thus, the work presented here can be viewed as a preliminary investigation into partial integro-differential equations of this type. We develop a first-order explicit finite difference method to solve the equation and is able to expose some interesting behavior of plasmons when the graphene's Drude weight is manipulated temporarily and spatially.

The body of this paper can be divided into three main parts: (i) Derivation of the current density equation, (ii) Finite difference method for solving the resulting equation, and (iii) Numerical experiments. In Section \ref{Section 2}, we describe the problem under consideration as an initial value problem for Maxwell's equations. Of interest is how a known initial field evolves as the Drude weight changes in time and space.

We derive a partial integro-differential equation for the current density on the conductive sheet in Section \ref{Section 3}. The Fourier transform is employed for this task. We consider the problem in the frequency domain and obtain a representation for the electric field that `lives' on the graphene sheet. After going back to the time domain we arrive at a reduced equation for the current density.

To solve the resulting partial integro-differential equation, we propose a finite difference scheme in Section \ref{Section 4}. Since we do not know what is the proper absorbing boundary conditions for the current density, we opt to solve the problem in a `padded' domain. We provide details of the grid construction. We show that the finite difference approximation is consistent, but defer the stability analysis to future work. The accuracy of the method is validated with known exact solutions. 

With the simulator at hand, we explore new phenomena associated with specific types of time- and space-dependent Drude weight. The paper ends with a discussion section.

\section{Modeling graphene with a conductive sheet} 
\label{Section 2}

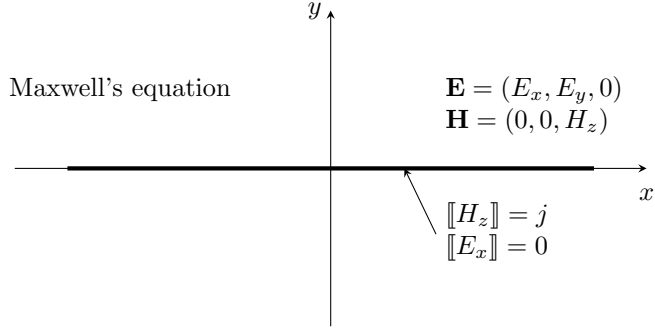
\begin{figure}
\begin{center}
\begin{tikzpicture}[scale=0.70]
\draw [-stealth] (0, 0) -- (12, 0);
\draw [-stealth] (6, -3) -- (6, 3);
\draw [ultra thick] (1, 0) -- (11, 0);
\draw (12,-0.5) node[]{$x$};
\draw (5.7, 3) node[]{$y$};
\draw (2., 1.5) node[]{Maxwell's equation};
\draw (8, 1.5) node[anchor=west]{${\bf E}=(E_x, E_y, 0)$};
\draw (8, 0.9) node[anchor=west]{${\bf H}=(0, 0, H_z)$};
\draw (8, -0.9) node[anchor=west]{$\llbracket H_z \rrbracket = j$};
\draw (8, -1.5) node[anchor=west]{$\llbracket E_x \rrbracket = 0$};
\draw (8, -1.25) [-stealth] -- (7.4, 0);
\end{tikzpicture}
\end{center}
\caption{Geometry of the Problem. The graphene sheet is placed along the $x$-axis. The TM fields are indicated. There are jump conditions on the sheet.}
\label{C2 Problem Setup}
\end{figure}

We consider a transverse magnetic (TM) electromagnetic field. With orthogonal coordinates $xyz$, the electric field is of the form $E = (E_x, E_y, 0)$ and the magnetic field is of the form $H = (0, 0, H_z)$, as considered in \cite{Maxwell_Fadil}. The fields are time and space dependent. We assume a constant permittivity $\epsilon$ and a constant permeability $\mu$. Figure \ref{C2 Problem Setup} provides a schematic of the problem geometry. The graphene sheet occupies the entire $x$-axis. The electromagnetic field away from the sheet is governed by Maxwell's equation in the TM mode. The sheet, being of zero thickness, is modeled by jump conditions on the $x$-axis. The goal of this work is to study how the field evolves when the Drude weight is varied temporally and spatially. 

For $y \ne 0$, we have the following simplified Maxwell's equation
\begin{align}
\mu \frac{\partial H_z}{\partial t} &= \frac{\partial E_x}{\partial y} - \frac{\partial E_y}{\partial x}, 
\label{Maxwell 1} \\
\epsilon \frac{\partial E_x}{\partial t} &= \frac{\partial H_z}{\partial y},
\label{Maxwell 2} \\
\epsilon \frac{\partial E_y}{\partial t} &= -\frac{\partial H_z}{\partial x}.
\label{Maxwell 3}
\end{align}

On the graphene sheet we have
\begin{equation}
\llbracket H_z \rrbracket_{y=0} = H_z(x, 0^+, t) - H_z(x, 0^-, t) = j(x, t),
\label{Hz jump}
\end{equation}
where $j(x,t)$ is the current density on the sheet. The tangential component of electric field is continuous across the graphene sheet
\begin{equation}
\llbracket E_x \rrbracket_{y=0} = 0.
\label{Ex jump}
\end{equation}
Drude's law provides a relationship between the current density and the $E_x$ field
\begin{equation}
\frac{\partial j}{\partial t} = -\frac{1}{\tau} j + D(x, t) E_x(x, 0, t) \label{Drude},
\end{equation}
where $\tau$ is the damping factor in the system and $D(x, t)$ is our time- and space-dependent Drude weight. 
From \eqref{Maxwell 2} and \eqref{Ex jump} we see that $\partial H_z/\partial y$ is continuous across $y=0$.

We analyze an initial value problem associated with the system. That is, we assume we are given $\{ E_x(x, y, 0), E_y(x, y, 0), H_z(x, y, 0) \}$ and hence $j(x, 0)$.  We will choose the initial data to be the fields associated with a constant Drude weight $D_0$. Existence and uniqueness results for initial value problems can be found \cite{Maxwell_Monk_Book}, for generic problem settings such as scattering in a heterogeneous medium and scattering by an object. We are not aware of such results for the problem at hand. Indeed, very little is known about wave phenomena when the material coefficients are time- and space-dependent. 

\subsection{Energy of the system}

We derive an energy balance equation and interpret the result in this section. A physics-based perspective on plasmon energy can be found in \cite{Energy_Bakunov}. Because of the discontinuity of $H_z$ along the $x$-axis, we consider the domain $\mathbb{R}^2$ as two half-spaces. We begin by multiplying \eqref{Maxwell 1} by $H_z$ and integrating over the upper half-space. After some integration-by-parts, we arrive at
\begin{align}
\frac{\mu}{2} \frac{d}{d t} \left[ \int_{-\infty}^\infty \int_{0^+}^\infty H_z^2 \,dydx \right] &= -\int_{-\infty}^\infty E_x(x, 0, t) H_z(x, 0^+, t) \,dx
\nonumber \\
&- \int_{-\infty}^\infty \int_{0^+}^\infty E_x \frac{\partial H_z}{\partial y} \,dydx - \int_{-\infty}^\infty \int_{0^+}^\infty \frac{\partial E_y}{\partial x} H_z \,dydx.
\nonumber
\end{align}
We perform a similar calculation for the lower half-space and add the resulting expression to the one above and get
\begin{align}
\frac{\mu}{2} \frac{d}{d t} \left[ \int_{\mathbb{R}^2} H_z^2 \,dydx \right] &= - \int_{-\infty}^\infty E_x(x, 0, t) \left[ H_z(x, 0^+, t) - H_z(x, 0^-, t) \right] \,dx
\nonumber \\
& \hspace{3em} - \int_{\mathbb{R}^2} E_x \frac{\partial H_z}{\partial y} \,dydx - \int_{\mathbb{R}^2} \frac{\partial E_y}{\partial x} H_z \,dydx \nonumber \\
& = - \int_{-\infty}^\infty E_x(x, 0, t) j(x, t) \,dx \nonumber \\
& \hspace{3em} - \int_{\mathbb{R}^2} E_x \frac{\partial H_z}{\partial y} \,dxdy + \int_{\mathbb{R}^2} E_y \frac{\partial H_z}{\partial x} \,dxdy
\nonumber
\end{align}
after applying \eqref{Hz jump} and performing integration-by-parts.

Multiplying \eqref{Maxwell 2} by $E_x$ and \eqref{Maxwell 3} by $E_y$, and integrating over $\mathbb{R}^2$, we get
\begin{equation*}
\frac{d}{dt} \int_{\mathbb{R}^2} \frac{1}{2}\epsilon \left( E_x^2 + E_y^2 \right) \,dxdy = \int_{\mathbb{R}^2} \frac{\partial H_z}{\partial y} E_x \,dxdy - \int_{\mathbb{R}^2} \frac{\partial H_z}{\partial x} E_y \,dxdy.
\label{Energy 2}
\end{equation*}
Summing we get the following equation
\[
\frac{d}{dt} \int_{\mathbb{R}^2} \left( \frac{1}{2}\mu H_z^2 + \frac{1}{2} \epsilon E_x^2 + \frac{1}{2} \epsilon E_y^2 \right) \,dxdy = - \int_{-\infty}^\infty E_x(x, 0, t) j(x, t) \,dx.
\]
Note that the bracketed term on the left is the energy density away from the graphene sheet. We now have to account for the energy stored in the graphene.

Solving $E_x(x, 0, t)$ from \eqref{Drude} and applying the result to the right-hand side, we have
\begin{align*}
- \int_{-\infty}^\infty & E_x(x, 0, t) j(x, t) \,dx =
 -\int_{-\infty}^\infty \frac{1}{D(x, t)} \frac{\partial j}{\partial t} j \,dx - \frac{1}{\tau} \int_{-\infty}^\infty \frac{1}{D(x, t)} j^2 \,dx \\
 &= -\frac{d}{dt} \int_{-\infty}^\infty \frac{1}{2D(x, t)} j^2 \, dx + \int_{-\infty}^\infty \left[ \frac{\partial}{\partial t} \frac{1}{2D(x, t)} \right] j^2 \,dx 
 - \frac{1}{\tau} \int_{-\infty}^\infty \frac{1}{D(x, t)} j^2 \,dx .
\end{align*}
Rearranging, we get
\begin{align*}
\frac{d}{dt} \int_{\mathcal{R}^2} \left( \frac{1}{2}\mu H_z^2 + \frac{1}{2} \epsilon E_x^2 + \frac{1}{2} \epsilon E_y^2 \right) & \,dxdy 
+ \frac{d}{dt} \int_{-\infty}^\infty \frac{1}{2D(x,t)}j^2 \, dx
= \\ & \int_{-\infty}^\infty \left[\frac{\partial}{\partial t} \frac{1}{2D(x, t)}\right] j^2 \,dx - \frac{1}{\tau} \int_{-\infty}^\infty \frac{1}{D(x, t)} j^2 \,dx.
\end{align*}
This equation can be viewed as an energy balance equation.  The first term on the left is the time rate of change of the energy in the vacuum away from the graphene, the second term is time rate of change of the energy stored in the graphene. Let us define energy in the system at time $t$ as
\begin{equation}
\mathcal{E}(t) =  \int_{\mathcal{R}^2} \left( \frac{1}{2}\mu H_z^2 + \frac{1}{2} \epsilon E_x^2 + \frac{1}{2} \epsilon E_y^2 \right) \,dxdy  
+ \int_{-\infty}^\infty \frac{1}{2D(x,t)}j^2 \, dx.
\label{Energy Definition}
\end{equation}
By integrating the energy balance equation above from $t=0$ to $t=T$, we arrive at
\begin{equation}
\mathcal{E}(T)-\mathcal{E}(0) = \frac{1}{2} \int_0^T \int_{-\infty}^\infty \left[ \frac{\partial}{\partial t} \left( \frac{1}{D(x, t)} \right) \right] j^2 \,dxdt
- \frac{1}{\tau} \int_0^T \int_{-\infty}^\infty \frac{1}{D(x, t)} j^2 \,dx.
\label{Energy Equation}
\end{equation}
On the right, we have the causes for the energy change. The first term on the right can be positive or negative depending on the time derivative of $1/D(x,t)$, whereas the second term is always negative and is caused by damping. In \eqref{Energy Equation}, if the Drude weight is time-independent, i.e., $D=D(x)$, and there is no damping, i.e., $\tau= \infty$, then energy is conserved. Varying $D$ in time can inject or remove energy in the system, and this was already shown in \cite{Josh_Physics_Review} \cite{Maxwell_Fadil}.

\subsection{Background solution}

We use the constant Drude weight solution from \cite{Maxwell_Fadil} as our background solution. Let $D_0$ be the Drude weight. The background solution is given by
\begin{align}
E_{x0}(x, y, t) &= \frac{\gamma_0}{\epsilon s_0} e^{i \xi_0 x} e^{-\gamma_0 |y|} e^{-s_0 t},
\label{Ex0 defi} \\
E_{y0}(x, y, t) &= \frac{i \xi_0}{\epsilon s_0} \sgn(y) e^{i \xi_0 x} e^{-\gamma_0 |y|} e^{-s_0 t},
\label{Ey0 defi} \\
H_{z0}(x, y, t) &= \sgn(y) e^{i \xi_0 x} e^{-\gamma_0 |y|} e^{-s_0 t},
\label{Hz0 defi} \\
j_0(x, t) &= 2 e^{i \xi_0 x} e^{-s_0 t}.
\label{j0 defi}
\end{align}
We view the horizontal wave number $\xi_0$ as the parameter. Then $\gamma_0$ and $s_0$ are determined from $\xi_0$. From the expressions above, we can see that $\gamma_0$ is the constant controlling the exponential decay of the field away from the graphene, whereas the real and imaginary parts of $s_0$ are the time-frequency and the damping constant of the field.

From \eqref{Maxwell 1} we have the dispersion relation
\begin{equation}
\gamma_0^2 = \mu \epsilon s_0^2 + \xi_0^2.
\label{gamma0 equa 1}
\end{equation}
Furthermore, Drude's law \eqref{Drude} reduces to a quartic equation for $s_0$
\begin{equation}
s_0^4 - \frac{2}{\tau} s_0^3 + \left( \frac{1}{\tau^2} - \frac{\eta^2 D_0^2}{4} \right) s_0^2 - \frac{D_0^2 \xi_0^2}{4 \epsilon^2} = 0,
\label{s0 equa}
\end{equation}
where $\eta = \sqrt{\mu / \epsilon}$ is the wave impedance. As shown in \cite{Maxwell_Fadil}, for large $\xi_0$, the desired $s_0$ with a nonzero imaginary part is given by
\begin{equation}
s_0 = i \sqrt{\frac{\xi_0 D_0}{2 \epsilon}} + \frac{1}{2 \tau} + O(\xi_0^{-\frac{1}{2}}).
\label{s0 approx}
\end{equation} 
If we set $\tau = \infty$, \eqref{s0 equa} simplifies to a quadratic equation for $s_0^2$. There are one positive and one negative root. The negative one yields a conjugate pair of purely imaginary roots for $s_0$. Those two real roots are not physically meaningful. 

\begin{figure}
\includegraphics[width=6.5cm]{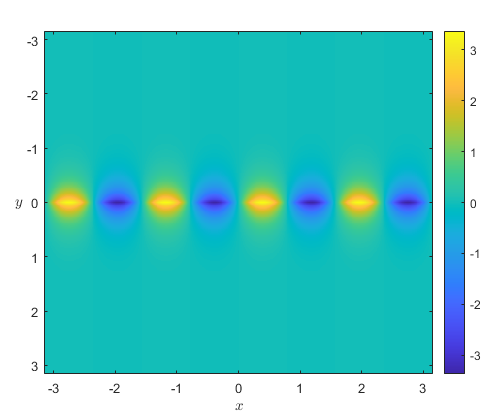}\includegraphics[width=6.5cm]{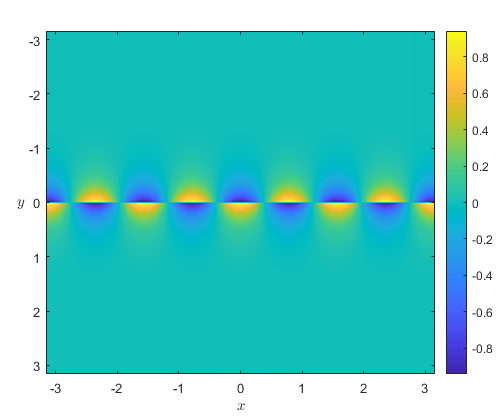}
\caption{Plots of the $E_{x0}$ (left) and $H_{x0}$ (right) for the constant Drude weight solution. Note that both fields are confined near the graphene sheet. The continuity in $E_x$ and the discontinuity in $H_z$ can be clearly seen.}
\label{C2 Background Ex0}
\end{figure}

We plot $E_{x0}(x, y, 0)$ and $H_{z0}(x, y, 0)$ on the box $[-3.15, 3.15] \times [-3.15, 3.15]$ in Figure \ref{C2 Background Ex0} as heat maps. For simplicity, we pick $\mu = 1$, $\epsilon = 1$, $\tau = \infty$, $D_0 = 0.675$ and $\xi_0 = 4$. We will justify this choice of parameters in Section \ref{Section 4}. We can see four full cycles, since the wavenumber is $\xi_0 = 4$ and the length of our box is close to $2 \pi$. We use the background solution at $\left\{E_{x0}(x,y,0), E_{y0}(x,y,0), H_{z0}(x,y,0), j_0(x,0)\right\}$ as initial values for \eqref{Maxwell 1}-\eqref{Drude}. This is an initial value problem driven by the Drude weight $D(x, t)$. We will focus on the evolution of the current density $j(x, t)$.

\section{An equation for the current density} 
\label{Section 3}

Analysis of the system when the Drude weight is time-dependent was greatly simplified when we were able to reduce the governing equations to a single equation on the conductive sheet \cite{Maxwell_Fadil}. Motivated by this observation, we next derive a new partial integro-differential equation for the current density $j(x, t)$ on the conductive sheet. The remainder of the paper focuses on this equation. 

Let $u(x,t)$ be some function which depends on time $t$ and space $x$. We will be using transforms methods -- Laplace transform in time and Fourier transform in space. We use the notation
\begin{equation*}
\widehat{u}(x,s) = \int_0^\infty u(x,t) e^{-s t} \,dt,
\end{equation*}
to denote the for time-Laplace transform of $u(x,t)$. By $\mathcal{L}_t$, we mean a Laplace transform operation in time.  The inverse operation is denoted by $\mathcal{L}_t^{-1}$.
The notation for Fourier transform in $x$ is
\begin{equation*}
  \widetilde{u}(\xi,t) =  \int_{-\infty}^\infty u(x,t) e^{-i \xi x} \,dx .
\end{equation*}
The Fourier transform operation and its inverse operation are denoted respectively by $\mathcal{F}_x$ and $\mathcal{F}_x^{-1}$. For example $\widetilde{\widehat{u}}(\xi,s) = \mathcal{F}_x \mathcal{L}_t \{u(x,t)\}$.

We turn to the derivation of an equation for $j(x, t)$ now. The core idea is to find a formula for $E_x(x, 0 , t)$ and insert into \eqref{Drude}. After taking the Laplace and Fourier transforms of the governing equations \eqref{Maxwell 2}-\eqref{Maxwell 3}, we solve for $\widetilde{\widehat{E_x}}(\xi,y,s)$ and $\widetilde{\widehat{E_y}}(\xi,y,s)$ which we substitute in \eqref{Maxwell 1} to get
\begin{equation}
\frac{\partial^2 \widetilde{ \widehat{H_z} }}{\partial y^2} (\xi, y, s) - \gamma^2 \widetilde{ \widehat{H_z} } (\xi, y, s) = f(\xi, y, s),
\label{LF Hz equa}
\end{equation}
where $\gamma^2 = \mu \epsilon s^2 + \xi^2$, and
\begin{equation}
    f(\xi, y, s) = -\mu \epsilon s \widetilde{H_z}(\xi, y, 0) - \epsilon \frac{\partial \widetilde{E_x}(\xi, y, 0)}{\partial y} 
+ i \xi \epsilon \widetilde{E_y}(\xi, y, 0).
\label{f defi}
\end{equation}
Following \cite{Maxwell_Fadil}, we solve \eqref{LF Hz equa} by multiplying the equation by $\mbox{sgn}(y) e^{-\gamma|y|}$ and integrating in $y$. The terms involving $\widetilde{\widehat{H_z}}$ can be integrated by part, producing 
\[
 2 \frac{\partial \widetilde{\widehat{H_z}}}{\partial y} (\xi,0,s) + \gamma \llbracket \widetilde{\widehat{H_z}} \rrbracket_{y=0} = -\int_{-\infty}^\infty f(\xi,y,s) \mbox{sgn}(y) e^{-\gamma|y|} \, dy.
\]
Recognizing that the jump in $H_z$ is the surface current, we arrive at
\begin{equation}
\frac{\partial \widetilde{ \widehat{H_z} }}{\partial y} (\xi, 0, s) = -\frac{\gamma}{2} \widetilde{ \widehat{j} } (\xi, s) - \frac{1}{2} \int_{-\infty}^\infty f(\xi, y, s) \sgn(y) e^{-\gamma |y|}  \,dy.
\label{dHzdy}
\end{equation}
Using the transformed version of \eqref{Maxwell 2} we have
\begin{align}
\widetilde{ \widehat{E_x} }(\xi, 0, s) &= \frac{1}{\epsilon s} \frac{\partial \widetilde{ \widehat{H_z} }}{\partial y}(\xi, 0, s) + \frac{1}{s} \widetilde{E_x}(\xi, 0, 0)
\nonumber \\
&= -\frac{\gamma}{2 \epsilon s} \widetilde{ \widehat{j} } (\xi, s) - \frac{1}{2 \epsilon s} \int_{-\infty}^\infty \sgn(y) e^{-\gamma |y|} f(\xi, y, s) \,dy + \frac{1}{s} \widetilde{E_x}(\xi, 0, 0).
\label{LF Ex}
\end{align}
The right hand side only involves $j(x, t)$ and known functions. We can rewrite \eqref{Drude}
\begin{equation}
\frac{\partial j}{\partial t} = -\frac{1}{\tau} j + D(x, t) \mathcal{F}_x^{-1} \mathcal{L}_t^{-1} 
\left\{ \widetilde{\widehat{E_x} } (\xi, 0, s) \right\}.
\label{Drude LF Ex}
\end{equation}
This will be an equation for the current density $j(x, t)$. The next section will focus on computing the inverse transforms.

\subsection{Inverse transforms}

To complete the equation for current density, we need to invert the transforms in \eqref{Drude LF Ex}. We have three terms in  \eqref{LF Ex}. For the first term we apply the convolution theorem for Fourier transforms. Recalling from \cite{Maxwell_Fadil}
\begin{equation}
\mathcal{L}_t^{-1} \left\{ \frac{\gamma}{s} \right\} = \frac{1}{c} k_1(\xi, t) + \frac{1}{c} \delta(t),
\label{gamma over s}
\end{equation}
where
\begin{equation}
k_1(\xi, t) = \int_0^t \frac{c \xi}{t'} J_1(c \xi t') \,dt'.
\label{k1 defi}
\end{equation}
Therefore, 
\begin{equation}
\mathcal{F}_x^{-1} \mathcal{L}_t^{-1} \left\{ -\frac{\gamma}{2 \epsilon s} \widetilde{ \widehat{j} } (\xi, s) \right\} 
= -\frac{\eta}{2} \mathcal{F}_x^{-1} \left\{ \widetilde{j} (\xi, t) + k_1(\xi, t) \,*_t\, \widetilde{j} (\xi, t) \right\}.
\label{Ex term 1}
\end{equation}
Here $c=1/\sqrt{\mu\epsilon}$ and $\eta=\sqrt{\mu/\epsilon}$.


To deal with the second term in \eqref{LF Ex}, we use the initial conditions, the fields at $t=0$ of \eqref{Ex0 defi}-\eqref{Hz0 defi}, calculate their $x$-Fourier transforms and insert them in \eqref{f defi} to get
\begin{equation}
f(\xi, y, s) = \left( -\mu \epsilon (s - s_0) - \frac{\xi_0}{s_0} (\xi - \xi_0) \right) 2 \pi \delta(\xi - \xi_0) e^{-\gamma_0 |y|} \sgn(y).
\label{f form}
\end{equation}
Note that $(\xi-\xi_0)\delta(\xi-\xi_0)$ is formally zero and $\int_{-\infty}^\infty e^{-(\gamma+\gamma_0)|y|}dy = -\frac{2}{(\gamma+\gamma_0)}$. So we have
\begin{align}
\mathcal{F}_x^{-1} \mathcal{L}_t^{-1} \left\{ - \frac{1}{2 \epsilon s} \int_{-\infty}^\infty \sgn(y) e^{-\gamma |y|} f(\xi, y, s) \,dy \right\} 
& = \mathcal{F}_x^{-1} \mathcal{L}_t^{-1} \left\{ \frac{\mu (s - s_0)}{s (\gamma + \gamma_0)} 2 \pi \delta(\xi - \xi_0) \right\}  .
\label{LF inverse}
\end{align}
 

Using the fact that $\mu\epsilon(s^2-s_0^2) = (\gamma^2 - \gamma_0^2) + (\xi^2 - \xi_0^2)$, we rewrite
\[
\frac{\mu (s - s_0)}{s (\gamma + \gamma_0)} = \frac{\mu (s-s_0)(s+s_0)}{s(\gamma+\gamma_0)(s+s_0)} = \frac{(\gamma - \gamma_0)}{\epsilon s (s + s_0)} - \frac{(\xi^2 - \xi_0^2)}{\epsilon s (s + s_0) (\gamma + \gamma_0)}.
\]
Notice that
\[
\mathcal{F}_x^{-1} \mathcal{L}_t^{-1} \left\{ \frac{\mu (s - s_0)}{s (\gamma + \gamma_0)} 2 \pi \delta(\xi - \xi_0) \right\} =
\mathcal{F}_x^{-1} \mathcal{L}_t^{-1} \left\{ \frac{(\gamma - \gamma_0)}{\epsilon s (s + s_0)} 2 \pi \delta(\xi - \xi_0) \right\},
\]
since the contribution from the second term is formally zero.

It can be shown that
\[
 \frac{(\gamma - \gamma_0)}{\epsilon s (s + s_0)} = \frac{\gamma_0}{\epsilon s_0} \left( \frac{1}{s + s_0}  -\frac{1}{s}\right) + \frac{1}{\epsilon} \frac{\gamma}{s} \frac{1} {(s + s_0)}.
\]
We can now easily invert the Laplace transform of the expression. The first terms are easy to invert while the second term requires using \eqref{gamma over s} and the convolution theorem. We get
\begin{align}
\mathcal{L}_t^{-1} \left\{ \frac{(\gamma - \gamma_0)}{\epsilon s (s + s_0)} \right\} &= \frac{\gamma_0}{\epsilon s_0} (e^{-s_0 t} - 1) + \frac{1}{\epsilon} \left(\frac{1}{c} k_1(\xi, t) + \frac{1}{c} \delta(t) \right) \,*_t\, e^{-s_0 t}
\nonumber \\
&= \frac{\gamma_0}{\epsilon s_0} (e^{-s_0 t} - 1) + \eta \left(k_1(\xi, t) \,*_t\, e^{-s_0 t} + e^{-s_0 t} \right).
\nonumber
\end{align}
We can take the inverse Fourier transform in \eqref{LF inverse}
\begin{align}
&\mathcal{F}_x^{-1} \mathcal{L}_t^{-1} \left\{ \frac{\mu (s - s_0)}{s (\gamma + \gamma_0)} 2 \pi \delta(\xi - \xi_0) \right\} \nonumber \\
&= \mathcal{F}_x^{-1} \left\{ \left[ \frac{\gamma_0}{\epsilon s_0} (e^{-s_0 t} - 1) + \eta \left(k_1(\xi, t) \,*_t\, e^{-s_0 t} + e^{-s_0 t} \right) \right] \delta(\xi-\xi_0) \right\} \nonumber \\
&= \frac{\gamma_0}{\epsilon s_0} e^{i \xi_0 x} (e^{-s_0 t} - 1) + \eta e^{i \xi_0 x} \left( k_1(\xi_0, t) \, *_t \, e^{-s_0 t} + e^{-s_0 t} \right).
\nonumber
\end{align}
Recalling the definition of $E_{x0}$ and $j_0$, we may identify and conclude
\begin{align}
\mathcal{F}_x^{-1} \mathcal{L}_t^{-1} & \left\{ \frac{\mu (s - s_0)}{s (\gamma + \gamma_0)} 2 \pi \delta(\xi - \xi_0) \right\} \nonumber\\
&= \frac{\eta}{2} (j_0 + k_1(\xi_0, t) \,\widehat{*}\, j_0) + E_{x0}(x, 0, t) - E_{x0}(x, 0, 0).
\label{Ex term 2}
\end{align}

The last term of \eqref{LF Ex} is straightforward, we have
\begin{equation}
\mathcal{F}_x^{-1} \mathcal{L}_t^{-1} \left\{ \frac{1}{s} \widetilde{E_x}(\xi, 0, 0) \right\} = E_x(x, 0, 0) = E_{x0}(x, 0, 0),
\label{Ex term 3}
\end{equation}
by the initial condition for $E_x$. 

We have now addressed the inverse transforms of all three terms in \eqref{LF Ex}.  The next step is to substitute what we have found, namely \eqref{Ex term 1}, \eqref{Ex term 2} and \eqref{Ex term 3}, in \eqref{Drude LF Ex}. After some rearrangements, we end up with
\begin{align}
\frac{\partial j}{\partial t} (x, t) &+ \left(\frac{1}{\tau} + \frac{\eta}{2} D(x, t)\right) j(x, t) + \frac{\eta}{2} D(x, t) \mathcal{F}_x^{-1} \left\{ k_1(\xi, t) \,*_t\, \widetilde{j} (\xi, t) \right\}
\nonumber \\
 &= \frac{\eta}{2} D(x, t) j_0(x, t) + \frac{\eta}{2} D(x, t) k_1(\xi_0, t)  \, *_t \, j_0(x, t) + D(x, t) E_{x0}(x, 0, t).
\label{Finv j equa}
\end{align}
The only task left is to evaluate the Fourier inverse on the left-hand side.

\begin{remark}
It is worth noting here that one can take the spatial Fourier Transform of \eqref{Finv j equa} and view the resulting equation as the governing equation for the Fourier transform of the current density $\tilde{j}(\xi,t)$. Using the convolution theorem, the resulting equation is
\begin{align*}
&\frac{\partial \widetilde{j}}{\partial t} (\xi, t) + \left(\frac{1}{\tau} + \frac{\eta}{2} \widetilde{D}(\cdot, t)\right) *_\xi \widetilde{j}(\cdot, t)(\xi) + \frac{\eta}{2} \widetilde{D}(\cdot, t) *_\xi \left\{ k_1(\cdot, t) \,*_t\, \widetilde{j} (\cdot, t) \right\}(\xi) =
\nonumber \\
 & \frac{\eta}{2} \widetilde{D}(\cdot, t) *_\xi \widetilde{j_0}(\cdot, t)(\xi) + \frac{\eta}{2} \widetilde{D}(\cdot, t)  *_\xi [k_1(\xi_0, t)  \, *_t \, \widetilde{j_0}(\cdot, t)](\xi) + \widetilde{D}(\cdot, t) *_\xi \widetilde{E_{x0}}(\cdot, 0, t)(\xi).
\end{align*}
Recalling that $j_0(x,t)$ and $E_{x0}(x,0,t)$ are known explicitly, we can simplify the right-hand side forcing terms. This equation bears some resemblance to the integro-differential integral equation for the amplitude of the plasmon in \cite{Maxwell_Fadil}. When the Drude weight is only a function of $t$, the equation simplifies. It is when in addition $\xi=\xi_0$ that we arrive at the simpler equation in \cite{Maxwell_Fadil}.
\end{remark}
\subsection{Partial integro-differential equation for current density}
In order to get an equation for $j(x,t)$, we need to evaluate the spatial Fourier transform on the left-hand side of \eqref{Finv j equa}. We write out the term
\begin{equation}
\mathcal{F}_x^{-1} \left\{ k_1(\xi, t) \,*_t\, \widetilde{j} (\xi, t) \right\} = \frac{1}{2 \pi} \int_{-\infty}^\infty e^{i \xi x} \int_0^t k_1(\xi, \tau) \Tilde{j}(\xi, t - \tau) \, d\tau d\xi .
\label{Finv layer 1}
\end{equation}
We interchange the order of integration and consider
\begin{equation*}
\int_{-\infty}^\infty e^{i \xi x} k_1(\xi, \tau) \Tilde{j}(\xi, t - \tau) \,d\xi  = \int_{-\infty}^\infty \int_0^\tau e^{i \xi x} \frac{c \xi}{t'} J_1(c \xi t') \Tilde{j}(\xi, t - \tau) \,dt' d\xi ,
\end{equation*}
after using \eqref{k1 defi}. Next, we apply the identity $J_1(z) = -J_0'(z)$ and change the order of integration to get
\begin{equation}
\int_{-\infty}^\infty e^{i \xi x} k_1(\xi, \tau) \Tilde{j}(\xi, t - \tau) \,d\xi = -\int_0^\tau \frac{1}{t'} \frac{d}{d t'} \int_{-\infty}^\infty e^{i \xi x} J_0(c \xi t') \Tilde{j}(\xi, t - \tau) \,d\xi dt'.
\label{Finv layer 2}
\end{equation}
The Fourier transform of Bessel function of zeroth order can be found in \cite{Integral_Transforms}. Using this we get
\[
\int_{-\infty}^\infty e^{i \xi x} J_0(c \xi t') \,d\xi = 2 (c^2 t'^2 - x^2)^{-\frac{1}{2}} I_{|x|< c t'},
\]
where $I_{|x|<ct'}$ is the indicator function for the set $\{|x|<ct'\}$.
Inserting this in \eqref{Finv layer 2}
\begin{align*}
& \int_{-\infty}^\infty e^{i \xi x} k_1(\xi, \tau) \Tilde{j}(\xi, t - \tau) \,d\xi \\ 
=& -2\int_0^\tau \frac{1}{t'} \frac{d}{d t'} \int_{-\infty}^\infty 
(c^2 t'^2 - x'^2)^{-\frac{1}{2}} I_{|x'|< c t'} \, j(x - x', t - \tau) \,dx'\, dt'\\ 
=& -2 \int_0^\tau \frac{1}{t'} \frac{d}{d t'} \int_{-c t'}^{c t'} (c^2 t'^2 - x'^2)^{-\frac{1}{2}} j(x - x', t - \tau) \,dx' \,dt' .
\end{align*}
To put the above expression into a useful form, we make a change of variable $x'=ct'\sin\theta$ and rewrite the integral in $d\theta$. After some manipulation, we arrive at
\begin{align*}
& \int_{-\infty}^\infty e^{i \xi x} k_1(\xi, \tau) \Tilde{j}(\xi, t - \tau) \,d\xi \\ 
=& 2 \int_0^\tau \frac{1}{t'} \int_{-c t'}^{c t'} \frac{x'}{t'} (c^2 t'^2 - x'^2)^{-\frac{1}{2}} \, j_x(x - x', t - \tau) \,dx' dt'.
\end{align*}

The inverse Fourier transform in \eqref{Finv layer 1} now becomes
\begin{align*}
\mathcal{F}_x^{-1} \left\{ k_1(\xi, t) \, *_t \, \widetilde{j} (\xi, t) \right\} &= \frac{1}{\pi} \int_0^t \int_0^\tau \int_{-c t'}^{c t'} \frac{x'}{t'^2} (c^2 t'^2 - x'^2)^{-\frac{1}{2}} \, j_x(x - x', t - \tau) \,dx' dt' d\tau \\
=& \frac{1}{\pi} \int_0^t \int_{-c t'}^{c t'} \frac{x'}{t'^2}  (c^2 t'^2 \!-\! x'^2)^{-\frac{1}{2}} \! \left[ \int_{t'}^t j_x(x - x', t - \tau) \,d\tau \right] dx' dt'   \\
=& \frac{1}{\pi} \int_0^t \int_{-c t'}^{c t'} \frac{x'}{t'^2} (c^2 t'^2\! - \!x'^2)^{-\frac{1}{2}} \! \left[ \int_{0}^{t - t'} \!\! j_x(x - x', \tau') \,d\tau' \right] dx' dt'.
\end{align*}
Defining
\begin{equation}
v(x, t) = \int_0^t j(x, \tau') \,d\tau',
\label{v defi}
\end{equation}
we have
\begin{align*}
\mathcal{F}_x^{-1} \left\{ k_1(\xi, t) \, *_t \, \widetilde{j} (\xi, t) \right\} &= \frac{1}{\pi} \int_0^t \int_{-c t'}^{c t'} \frac{x'}{t'^2}  (c^2 t'^2 - x'^2)^{-\frac{1}{2}} \, v_x(x-x',t-t') dx' dt'\\
&= - \frac{1}{\pi} \int_0^t \frac{1}{t'^2} \int_{-c t'}^{c t'} (c^2 t'^2 - x'^2)^{\frac{1}{2}} \, v_{xx}(x - x', t - t') \,dx' dt',
\end{align*}
after an integration-by-parts. 

To simplify notations, for a function $g(x,t)$ we define the operator
\begin{equation}
\Phi[g] (x, t) = \int_0^t \frac{1}{t'^2} \int_{-c t'}^{c t'} ( c^2 t'^2 - x'^2)^{\frac{1}{2}} \, g_{xx}(x - x', t - t') \,dx' dt'.
\label{Phi defi}
\end{equation}
The equation for $v(x,t)$ is obtained from \eqref{Finv j equa} and takes the form of a partial integro-differential equation
\begin{align}
\frac{\partial^2 v}{\partial t^2} (x, t) &+ \left( \frac{1}{\tau} + \frac{\eta}{2} D(x, t) \right) \frac{\partial v}{\partial t} (x, t) - \frac{\eta}{2 \pi} D(x, t) \Phi[v] (x, t)
\nonumber \\
&= \frac{\eta}{2} D(x, t) j_0(x, t) + \frac{\eta}{2} D(x, t) k_1(\xi_0, t) \,*_t\, j_0(x, t) + D(x, t) E_{x0}(x, 0, t).
\label{v equa}
\end{align}
Our main concern is solving the initial value problem associated with \eqref{v equa} where $D=D_0$ for $t=0$, and $v(x,0)=0$, and $v_t(x,0)=j_0(x,0)$. It can be seen that \eqref{v equa} is nonlocal in time and space. The value of $\Phi[v](x,t)$ depends on $v_{xx}$ behind the backwards light cone from $(x,t)$ as \eqref{Phi defi} shows. This situation is illustrated in Figure \ref{nonlocal}. The kernel in $\Phi[\cdot]$ is finite for $t'\neq 0$, and its contribution is smallest for $t'=t$. Therefore, the convolution tells us that the most recent values of $v_{xx}$ are the most important.
\begin{figure}
\begin{center}
\begin{tikzpicture}[scale=0.85]
\draw [->] (-2,0) -- (10,0);
\draw [->] (0,0) -- (0,5);
\draw (-0.4, 5) node[]{$t$};
\draw (10.4, 0) node[]{$x$};
\draw (-1,0) -- (4,4);
\draw (4,4) -- (9,0);
\fill (4,4) circle (1.5pt);
\draw (4.5,4.2) node[]{$(x,t)$};
\draw (-1,-0.4) node[]{$x-ct$};
\draw (9,-0.4) node[]{$x+ct$};
\draw [dashed] (0.9, 1.5) -- (7.1, 1.5);
\draw [dashed] (0.9, 1.5) -- (0.9, 0);
\draw [dashed] (7.1, 1.5) -- (7.1, 0);
\draw (0.9, -0.4) node[]{$x-ct'$};
\draw (7.1, -0.4) node[]{$x+ct'$};
\draw (0, 1.5) -- (0.2, 1.5);
\draw (-0.4, 1.5) node[]{$t'$};

\end{tikzpicture}
\end{center}
\caption{Nonlocal property of the operator in \eqref{v equa}. The double integral in the operator is over the backwards light cone at $(x,t)$. Integration in $x'$ is over the horizontal dashed lines, and $t'$ is from $0$ to $t$.}
\label{nonlocal}
\end{figure}
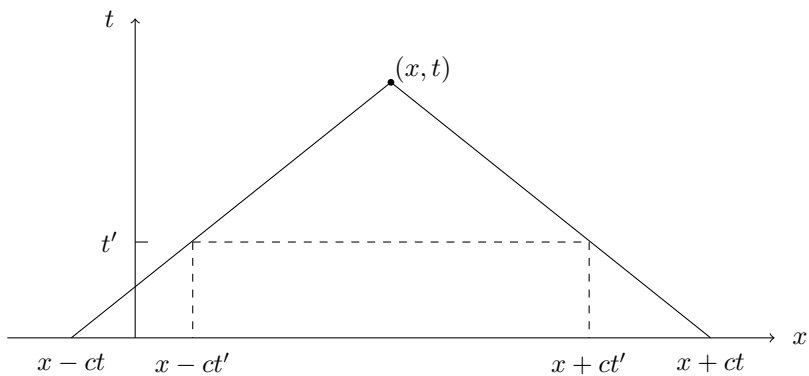
\begin{remark}
It is not obvious how one should go about establishing well-posed-ness of this initial value problem. We can verify that if $D(x,t)=D_0$, the solution is $v(x,t)=\int_0^t j_0(x,\tau')d\tau'$. We make two observations on this time-dependent integro-differential equations. First, it is not apparent from its form that this has a wave-like solution. Second, the Drude weight appears in the second term as part of the damping. This is somewhat counter-intuitive since if Drude weight is constant and $\tau = \infty$, there is no damping. 
\end{remark}

The remainder of this paper is focused on developing a finite difference method to solve the initial value problem for $v(x,t)$. We will also provide numerical simulations of the current density $j(x,t)$ for a specific form of $D(x,t)$. These simulations exhibit interesting interactions between the Drude weight and the current density.
\section{Finite difference method for the reduced equation} 
\label{Section 4}

In this section, we introduce a consistent finite difference method to solve the reduced equation for $v(x, t)$, given the initial data and $D(x, t)$. Of interest is the current density $j(x,t)$ which is the time derivative of $v(x,t)$.

We will work with a dimensionless version of \eqref{v equa}. We will use characteristic length $L=1$ \emph{micron}, and characteristic time $T = L/c =3.33$ \emph{femtosecond} ($10^{-15}$ \emph{s}). Under this rescaling, the damping factor is now $\tau/T$.  If we scale $E$ by $\overline{E}$ and $H$ and $j$ by $\overline{H}$, with $\overline{E}/\overline{H}=\sqrt{\mu/\epsilon}$, then the only parameter left in the problem is a dimensionless version of $D(x,t)$, which is $D$ multiplied by $T\sqrt{\mu/\epsilon}$. With $D_0=5.373 \times 10^{11}$ \emph{F}/\emph{s$^2$}, the rescaled $D_0 = 0.675$. The rescaled system is now
\begin{align}
\frac{\partial^2 v}{\partial t^2} (x, t) &+ \left( \frac{1}{\tau} + \frac{1}{2} D(x, t) \right) \frac{\partial v}{\partial t} (x, t) - \frac{1}{2 \pi} D(x, t) \Phi[v] (x, t)
\nonumber \\
&= \frac{1}{2} D(x, t) j_0(x, t) + \frac{1}{2} D(x, t) k_1(\xi_0, t) \,*_t\, j_0(x, t) + D(x, t) E_{x0}(x, 0, t).
\label{v equanondim}
\end{align}
Here, the nondimensional version of the operator $\Phi$ is
\begin{equation}
\Phi[g] (x, t) = \int_0^t \frac{1}{t'^2} \int_{-t'}^{t'} ( t'^2 - x'^2)^{\frac{1}{2}} \, g_{xx}(x - x', t - t') \,dx' dt'.
\label{Phi definondim}
\end{equation}
Also, the kernel $k_1$ is now
\[
k_1(\xi,t) = \int_0^t \frac{\xi}{t'}J_1(\xi t') dt'.
\]
We will be using $\xi_0=4$ which corresponds to a spatial wavelength of $\pi/2$ \emph{microns}.
\subsection{Finite difference method}

\begin{figure}
\begin{center}
\begin{tikzpicture}[scale=0.85]
\draw (0, 0) -- (12, 0);
\draw (3, 0) -- (3, 3);
\draw (9, 0) -- (9, 3);
\draw (3, 3) -- (9, 3);
\draw (0, 0) -- (3, 3);
\draw (9, 3) -- (12, 0);
\draw [->] (6, 0) -- (6, 3.5);
\draw [->] (12, 0) -- (13, 0);
\draw (6.3, 3.5) node[]{$t$};
\draw (12.9, -0.4) node[]{$x$};

\draw [dashed]  (0, 0) -- (0, -0.7);
\draw [dashed]  (12, 0) -- (12, -0.7);
\draw [dashed]  (3, 3) -- (0, 3);
\draw (-0.5, 3) node[]{$N \Delta t$};
\draw (0, -1) node[]{$-(M_1 + N) \Delta x$};
\draw (3, -0.5) node[]{$-M_1 \Delta x$};
\draw (6, -0.5) node[]{$0$};
\draw (9.18, -0.5) node[]{$M_1 \Delta x$};
\draw (11.7, -1) node[]{$(M_1 + N) \Delta x$};

\foreach \j in {0,1,...,6} {
    \foreach \i in {6,...,18} {
        \fill (\i/2,\j/2) circle (1pt);
    }
}
\foreach \n in {1,2,...,6} {
        \foreach \m in {1,...,\n} {
            \fill (\m/2-\n/2+2.5,\m/2-1/2) circle (1pt);
    }
}
\foreach \n in {1,2,...,6} {
        \foreach \m in {1,...,\n} {
            \fill (\n/2-\m/2+9.5,\m/2-1/2) circle (1pt);
    }
}

\draw (5, 0) -- (7.5, 2.5);
\draw (7.5, 2.5) -- (10, 0);
\fill (7.5, 2.5) circle (1.5pt);
\draw [dashed] (7.5, 2.5) -- (0, 2.5);
\draw (-0.5, 2.5) node[]{$k \Delta t$};
\draw [dashed] (7.5, 0) -- (7.5, 2.5);
\draw (7.5, -0.5) node[]{$l \Delta x$};

\draw [dashed] (8.5, 1.5) -- (1, 1.5);
\draw (-0.2, 1.5) node[]{$(k - k')\Delta t$};
\fill (7.5, 1.5) circle (1.5pt);
\fill (6.5, 1.5) circle (1.5pt);
\fill (8.5, 1.5) circle (1.5pt);

\draw [dashed] (6.5, 1.5) -- (6.5, -1.4);
\draw (5.5, -1.7) node[]{$l \Delta x - k' \Delta x$};
\draw [dashed] (8.5, 1.5) -- (8.5, -1.4);
\draw (9.5, -1.7) node[]{$l \Delta x + k' \Delta x$};

    
\end{tikzpicture}
\end{center}
\caption{Illustration of the grid system used in the computation. Of interest is the $v(x, t)$ for $|x| \leq A = M_1 \Delta x$. After discretization, this means we wish to calculate $v_{kl}$ for $l\in[-M_1 : M_1]$. Due to the dependence of the solution $v(x, t)$ on its past, we need a trapezoidal region as indicated. For a pair $(k, l)$, the dependence of the solution at this point is all the solutions in the triangular region indicated in the figure.}
\label{C4 FDM Grid}
\end{figure}
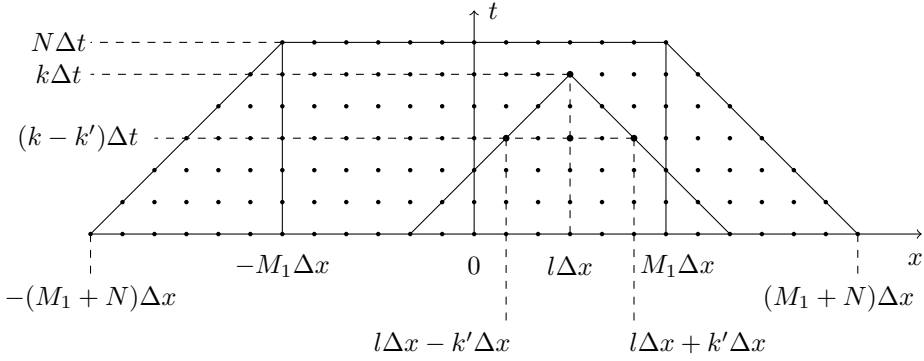

Our goal is to solve an initial value problem associated with \eqref{v equanondim}. We are given $D(x, t)$ and the initial data on the whole real line. The domain of interest is $(x,t)$ in $[-A,A] \times [0, T]$. The rectangular domain is regularly meshed with spacings of $\Delta x$ and $\Delta t$. We use indices $(l,k)$ to locate a grid point in the $(x,t)$ space. We let
\[
v_{l, k} \approx v(l\Delta x, k\Delta t).
\]
Because of the dependence of $v$ at $(x=l\Delta x,t=k\Delta t)$ on its values behind the backwards light cone, to simplify the computation of the operation $\Phi[v](x,t)$, we set the ratio $\Delta t/\Delta x$ to be an integer. Since we are modeling a wave phenomenon, it is unlikely that we will achieve a desired level of accuracy for $\Delta t/\Delta x$ greater than 1. So, for the remainder of this work, we will set it to 1. This choice makes sure that the light rays emanating from $(x,t)$, i.e., lines of slope $\pm 1$ at $(x,t)$ are on grid points. Moreover, with this choice, a horizontal line $t=k\Delta t$ contained inside the light cone has an integer number of grid points (see Figure \ref{C4 FDM Grid}). 

We set $A=M_1 \Delta x$ and $T=N\Delta t$.  We have one more detail to attend to. We want to avoid prescribing boundary conditions at $|x|=A$. We do this by padding the domain of interest.  Since the solution at $(\pm A,T)$ depend on values of $v$ in the backwards light cone, we pad the domain with triangular regions as shown in Figure \ref{C4 FDM Grid}. The resulting computational domain is a trapezoid, with base $[-(M_1+N)\Delta x, (M_1+N)\Delta x ]$, top $[-M_1 \Delta x, M_1\Delta x]$, and height $N\Delta t$.  



Next we devise an explicit time marching scheme. Let $x=l\Delta x$ and $t=k\Delta t$; i.e., we are at grid point $(l,k)$ and the quantity of interest is $v_{l,k}$. We denote right-hand side of \eqref{v equanondim} at this grid point by $R_{l, k}$. 
For the left-hand side, we define
\begin{align}
A_{l, k} &= \frac{1}{\tau} + \frac{1}{2} D(l \Delta x, k \Delta t),
\nonumber \\
B_{l, k} &= -\frac{1}{2 \pi} D(l \Delta x, k \Delta t).
\nonumber
\end{align}
We use $F_{l, k}$ to denote the approximate value of $\Phi[v](l\Delta x,k\Delta t)$.
For partial derivatives, we apply the central difference approximations.
Then \eqref{v equanondim} is approximated by
\begin{equation}
\frac{v_{l,k + 1} + v_{l,k - 1} - 2v_{l,k}}{\Delta t^2} + A_{l,k} \frac{v_{l,k + 1} - v_{l,k - 1}}{2 \Delta t} + B_{l,k} F_{l,k} = R_{l,k}.
\label{diff equa}
\end{equation}
We further define
\begin{align}
A_{l,k}^+ &= \frac{A_{l,k}}{2 \Delta t} + \frac{1}{\Delta t^2},
\nonumber \\
A_{l,k}^- &= \frac{A_{l,k}}{2 \Delta t} - \frac{1}{\Delta t^2}.
\nonumber
\end{align}
From \eqref{diff equa} we obtain an explicit time-marching scheme
\begin{equation}
v_{l,k + 1} = \frac{1}{A_{l,k}^+} \left( R_{l,k} - B_{l,k} F_{l,k} + \frac{2}{\Delta t^2} v_{l,k} + A_{l,k}^- v_{l,k - 1} \right).
\label{update form}
\end{equation}
At first glance, the time-marching scheme appears to involve the node values of the most recent two time layers. However, we should note that $F_{l,k}$ depends on values of $v_{l,k}$ behind the light cone at $(x_l,t_k)$. We will treat the calculation of $F_{l,k}$ momentarily.

To initialize the scheme, we need $v_{l,0}$ and $v_{l,1}$.
It is easy to see that $v_{l,0} = 0$ and $F_{l,0} = 0$ for all $l$. Therefore, at time layer $k = 0$, \eqref{update form} reads
\begin{equation}
v_{l,1} = \frac{1}{A_{l,0}^+} \left( R_{l,0} + A_{l,0}^- v_{l,-1} \right).
\label{first layer 1}
\end{equation}
We know that $v_t(x, 0) = j_0(x, 0)$ so it is known. We can evaluate 
\[
Q_l = j_0(l\Delta x, 0).
\]
Using the central difference approximation we have
\begin{equation}
\frac{v_{l,1} - v_{l,-1}}{2 \Delta t} = Q_l.
\label{first layer 2}
\end{equation}
Equations \eqref{first layer 1} and \eqref{first layer 2} form a linear system for $v_{l,1}$ and $v_{l,-1}$. We solve it to get
\begin{equation}
v_{l,1} = \frac{1}{A_{l,0}^+ - A_{l,0}^-} \left( R_{l,0} - 2 \Delta t A_{l,0}^- Q_l \right) = \frac{\Delta t^2}{2} \left( R_{l,0} - 2 \Delta t A_{l,0}^- Q_l\right).
\label{first layer}
\end{equation}
With this, we now have $v_{l,0}$ and $v_{l,1}$ and can march forward in time.

\subsection{Approximation of the integral operator}

We focus on the $\Phi[v](x,t)$ term in \eqref{v equanondim}, defined in \eqref{Phi definondim}. We use
\[
F_{l,k} \approx \Phi[v](l\Delta x,k \Delta t).
\]
From \eqref{Phi definondim} we see that $F_{k,l}$ depends on $v_{xx}(x,t)$ for all $(x,t)$ behind the light cone at $(l\Delta x, k\Delta t)$. By our choice of $\Delta t = \Delta x$, this dependence can be handled easily. We use
\[
t' = k' \Delta t, \;\; x' = l' \Delta x.
\]
We introduce the notation
\[
G(l,k,k') \approx \frac{1}{k'^2 \Delta t^2} \int_{- k' \Delta t}^{ k' \Delta t} (k'^2 \Delta t^2 - x'^2)^{\frac{1}{2}} v_{xx}(l \Delta x - x', k \Delta t - k' \Delta t) \,dx'.
\]
We replace the outer integral in \eqref{Phi definondim} by the trapezoidal rule
\begin{equation*}
    F_{l,k} 
= \Delta t \left[ \frac{G(l,k, 0)}{2} + \frac{G(l,k, k)}{2} + \sum_{k' = 1}^{k - 1} G(l,k, k') \right].
\end{equation*}

Next, we let
\begin{equation}
P(l, l', k, k') = (k'^2 \Delta t^2 - l'^2 \Delta x^2)^{\frac{1}{2}} v_{xx}(l \Delta x - l' \Delta x, k \Delta t - k' \Delta t).
\label{Peqn}
\end{equation}
Then, using trapezoidal rule, we have
\begin{align*}
&G(l,k, k')  \\
&= \frac{\Delta x}{k'^2 \Delta t^2} \left[ \frac{P ( l, k', k, k', ) }{2} + \frac{P ( l, -k', k, k') }{2} + \sum_{l' = - k' + 1}^{k' - 1} P(l, l',k, k') \right].
\end{align*}
We cannot use this equation to determine $G(l,k, 0)$. Instead, we use 
\begin{align*}
    G(l, k, 0) &= 
\lim_{t' \to 0} \frac{1}{t'^2} \int_{- t'}^{ t'} ( t'^2 - x'^2)^\frac{1}{2} v_{xx}(l \Delta x - x', k \Delta t - t') \,dx' \\
&= \lim_{t' \to 0} \int_{-1}^{1} (1 - \bar{x}^2)^\frac{1}{2} v_{xx}(l \Delta x -  t' \bar{x}, k \Delta t - t') \,d\bar{x} \\
&= \frac{1}{2} \pi v_{xx}(l \Delta x, k \Delta t).
\end{align*}
Notice in \eqref{Peqn} that
\[
P(l,k',k,k') = P(l,-k', k,k')=0.
\]
Therefore, the formula for $G(l,k, k')$ with $k' \neq 0$ simplifies to 
\[
G(l,k, k') = \frac{\Delta x}{k'^2 \Delta t^2} \sum_{l' = - k' + 1}^{k' - 1} P(l,l',k, k').
\]
Combining the results above, we can get a more compact formula for $F_{l,k}$
\begin{align}
F_{l,k} &= \frac{\pi c^2 \Delta t}{4} v_{xx}(l\Delta x, k\Delta t) + \frac{1}{2 k^2} \Delta x \sum_{l' = - k + 1}^{k - 1} (k^2 - l'^2 )^{\frac{1}{2}} v_{xx}((l - l')\Delta x,0)
\nonumber \\
&+ \Delta x \sum_{k' = 1}^{k - 1} \frac{1}{k'^2} \sum_{l' = - k'\frac{M_1}{N} + 1}^{k' - 1} (k'^2 - l'^2 )^{\frac{1}{2}} v_{xx}( (l - l')\Delta x, (k - k')\Delta t).
\nonumber
\end{align}
Since $v(x,0)=0$, we should have $v_{xx}(l\Delta x, 0) = 0$ for all $l$. This further simplifies the formula to
\begin{equation}
F_{l,k} = \frac{\pi \Delta t}{4} v_{xx}(l\Delta x, k\Delta t) + \Delta x \sum_{k' = 1}^{k - 1} \frac{1}{k'^2} \sum_{l' = - k' + 1}^{k' - 1} (k'^2 - l'^2 )^{\frac{1}{2}} v_{xx}( (l - l')\Delta x, (k - k')\Delta t ).
\label{Fkj form}
\end{equation}
The partial derivative $v_{xx}(l\Delta x,k\Delta t)$ can be approximated by finite differences from the values of $v$ at grid points. The form of \eqref{Fkj form} requires storage of previously computed $v_{xx}(l\Delta x,k\Delta t)$, which may pose an issue for large problems. This completes the derivation of a formula to approximate the operator $\Phi[v]$. 

\subsection{Consistency}

We wish to establish that the finite difference method in \eqref{update form} with \eqref{Fkj form} is consistent with the integro-differential equation \eqref{v equanondim} and \eqref{Phi definondim}. To this end, let $v(x,t)$ be a smooth function. We use the notation $d^+_t$ and $d^-_t$ to denote forward and backward difference approximation of time derivative, and $d^0_t$ to denote the central difference approximation of the same.  The truncation error we wish to estimate is
\begin{align}
{\mathcal E}(x_l,t_k) = \frac{\partial^2 v}{\partial t^2} (x_l, t_k) + &\left( \frac{1}{\tau} + \frac{1}{2} D(x_l, t_k) \right) \frac{\partial v}{\partial t} (x_j, t_k) - \frac{1}{2 \pi} D(x_l, t_k) \Phi[v] (x_l, t_k) \nonumber \\
&- (d^+_t d^-_t) v(x_l,t_k) - A_{l,k} d^0_t \, v(x_l,t_k) 
-  B_{l,k} F_{l,k}.
\label{consis error}
\end{align}
Using Taylor series, we know that
\begin{align*}
(d^+_t d^-_t) v(x_l,t_k) &= \frac{\partial^2 v}{\partial t^2}(x_l,t_k) + O(\Delta t^2),\\
d^0_t v(x_l,t_k) &= \frac{\partial v}{\partial t} (x_l,t_k) + O(\Delta t^2).
\end{align*}
So, from \eqref{consis error} we have
\begin{equation*}
{\mathcal E}(x_l,t_k) = O(\Delta t^2) - \frac{1}{2 \pi} D(x_l, t_k) \Phi[v] (x_l, t_k) - B_{l,k} F_{l,k}.
\end{equation*}
Equation \eqref{Fkj form} is a trapezoidal rule approximation of the double integral in \eqref{Phi definondim}. Therefore, we commit a truncation error of $O(\Delta x^2 + \Delta t^2)$. However, we need to further replace $v_{xx}(x,t)$ in \eqref{Fkj form} with its finite difference approximation, which produces an error bounded by $M \Delta x^2$, for some $M$. Let us write
\begin{align*}
  & \left| -\frac{1}{2 \pi} D(x_l, t_k) \Phi[v] (x_l, t_k) - B_{l,k} F_{l,k} \right|\\
& \leq K(\Delta x^2 + \Delta t^2) + \frac{\pi c^2 \Delta t}{4} M\Delta x^2 + \Delta x \sum_{k' = 1}^{k - 1} \frac{1}{k'^2} \sum_{l' = - k' + 1}^{k' - 1} (k'^2 - l'^2 )^{\frac{1}{2}} M \Delta x^2\\
&  \leq K(\Delta x^2 + \Delta t^2) + \frac{\pi c^2 \Delta t}{4} M\Delta x^2 + M \Delta x^3\sum_{k' = 1}^{k - 1} \frac{1}{k'^2} \sum_{l' = - k' + 1}^{k' - 1} |k'|.
\end{align*}
The first term represents the quadrature error for some $K > 0$. Note that 
\[
M \Delta x^3\sum_{k' = 1}^{k - 1} \frac{1}{k'^2} \sum_{l' = - k' + 1}^{k' - 1} |k'| \leq M \Delta x^3\sum_{k' = 1}^{k - 1} \frac{1}{k'^2} |k'| (2 k') \leq  M \Delta x^3 k.
\]
Therefore, we have
\begin{align*}
&\left| - \frac{1}{2 \pi} D(x_j, t_k) \Phi[v] (x_l, t_k) - B_{l,k} F_{l,k} \right| \\
& \hspace{1cm} \leq K(\Delta x^2 + \Delta t^2) + \frac{\pi c^2 \Delta t}{4} M\Delta x^2 + 2 M \Delta x^2 T,
\end{align*}
since $N\Delta t=T$ and $k\leq N$. Therefore, we have
\begin{equation}
  {\mathcal E}(x_l,t_k) = O(\Delta t^2) + O(\Delta x^2),  
\end{equation}
and the finite difference approximation is consistent.

\subsection{A Numerical Study of Convergence}
We compare the analytical solution for the current density \eqref{j0 defi} when the Drude weight is constant with numerically calculated approximate solution for different $\Delta x$. We consider the real part of the current density
\[
\mathcal{R}[j_0](x, t) = 2 \cos{(\xi_0 x - \omega_0 t)}.
\]
To produce the finite difference approximation of the solution, we use $j(x,0)=2\cos \xi_0 x$ and \eqref{update form} with \eqref{Fkj form}. In the calculations, we use $D(x,t)=D_0=0.675$ and consider the zero damping case $\tau=\infty$. The wavelength parameter $\xi_0=4$ and compute $\omega_0$ using \eqref{s0 equa}. With these parameter values $\omega_0 = 1.1377$. We compare the solution for the current density with the exact solution at time $T$. 

\begin{figure}
\centering
\includegraphics[width=6cm]{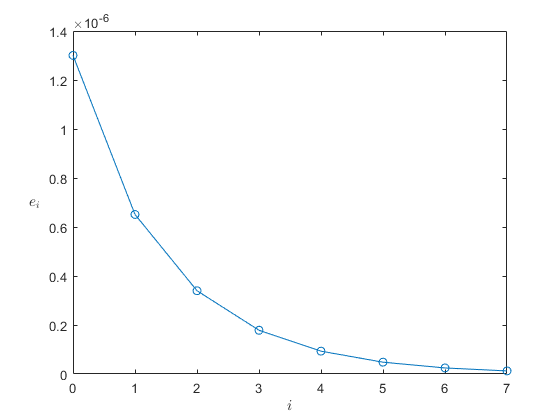}
\includegraphics[width=6cm]{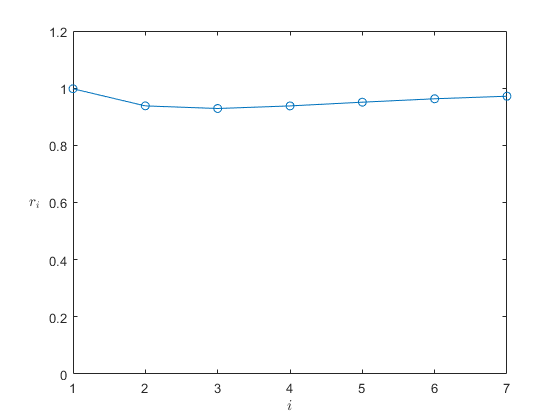}
\caption{The left figure shows the absolute error $e_i$ as a function of $i$ which controls the mesh size $\Delta x = 0.0105 \times 2^{-i}$. The right figure shows the estimated order of convergence as a function of $i$. }
\label{C4 Convergence Error Curve}
\end{figure}

We use obtain the finite difference solution at different mesh sizes. For the $i^{th}$ experiment, we set $\Delta x = 0.0105 \times 2^{-i}$, $M_1 = 5 \times 2^{i}$. We keep $\Delta t = \Delta x$ and set $N = 10 \times 2^{i}$. This means the grid is halved in the next experiment, while keeping the $A = M_1 \Delta x = 0.0525$ and $T = N \Delta t = 0.105$ unchanged. For each experiment, we record the $L^2$ error at $T = N \Delta t$ 
\begin{align}
e^2_i &= |v_{N, -M_1} - v(-A, T)|^2 \cdot \frac{\Delta x}{2} + |v_{N, M_1} - v(A, T)|^2 \cdot \frac{\Delta x}{2}
\nonumber \\
&+ \sum_{j = -M_1 + 1}^{M_1 - 1} |v_{N, j} - v(j \Delta x, T)
|^2 \cdot \Delta x,
\label{C4 convergence error}
\end{align}
and define the estimated order for $i \ge 1$
\begin{equation}
r_i = \frac{\log(e_{i - 1} / e_i)}{\log(2)}.
\label{C4 convergence order}
\end{equation}
The results are exhibited in Figure \ref{C4 Convergence Error Curve} where we show both the absolute $L^2(-A,A)$ and the estimated order of convergence.

From these figures, we conclude that the error of our scheme is $O(\Delta x)$, or the order of accuracy is 1.

\section{Numerical Experiments} 
\label{Section 5}
We conducted a number of numerical experiments with specific forms of the Drude weight $D(x,t)$. For time-dependent Drude weights, we were able to reproduce the phenomena first observed in \cite{Josh_Physics_Review}. In the following, we consider Drude weight of the form
\begin{equation}
D(x,t) = D_0 + \alpha \cos(\xi_1 x - \omega_1 t).
    \label{wavyDrude}
\end{equation}
The initial current density is
\[
j(x,0) = 2 \cos \xi_0 x.
\]
If the Drude weight is constant for $t>0$, the corresponding solution would be $j(x,t)=2\cos(\xi_0 x - \omega_0 t)$. Instead of being a constant, the Drude weight is $D_0$ plus a perturbation. The perturbation has a form similar to the background current density, i.e., it is wave-like. We have control over the perturbation's wave number $\xi_1$ and its frequency $\omega_1$.

In our experiments, we fix $\alpha=0.02$, $\xi_0=4$ and $D_0=0.675$ with no damping. This leads to $\omega_0=1.1377$. We will display $j(x,t)-j_0(x,t)$, which corresponds to the perturbation in the current density from the background.

\begin{figure}
\centering
\includegraphics[width=8cm]{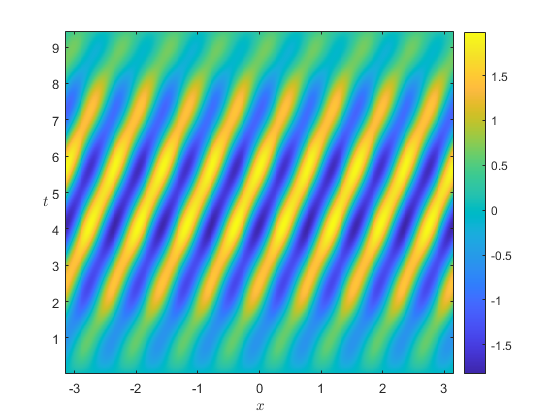}
\caption{Plot of the perturbational current density. Note that we have created a traveling perturbational current density.}
\label{FDM-moving}
\end{figure}

In the first experiment, we choose $\xi_1=\xi_0$ and $\omega_1=\omega_0$. This means the perturbation in the Drude weight is a wave with the same wave number and frequency and traveling in the same direction as the background current density. The resulting perturbational current density is a wave traveling to the right as displayed in Figure \ref{FDM-moving}. The phase speed is controlled by the choice of parameters $\omega_1$ and $\xi_1$.

\begin{figure}
\centering
\includegraphics[width=8cm]{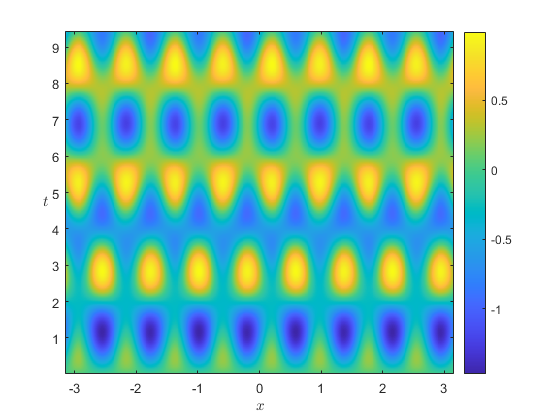}
\caption{Plot of the perturbational current density. Here, we have created a perturbational current density that is standing.}
\label{FDM-standing}
\end{figure}
In the second experiment, our goal is to create a standing current density perturbation. To do this, we pick $\xi_1=\xi_0$ and $\omega_1=-\omega_0$. This means that the Drude weight perturbation is a wave with the same wave number and frequency as the background current density but traveling in the opposite direction. The resulting perturbational current density is shown in Figure \ref{FDM-standing}. 

\begin{figure}
\centering
\includegraphics[width=8cm]{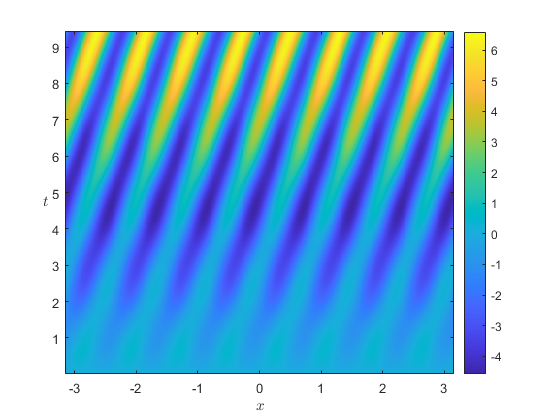}
\caption{Plot of the perturbational current density. With the parameters chosen from an approximate analysis in \cite{tong-thesis} we have created a perturbational current density that is both traveling and growing.}
\label{FDM-growing}
\end{figure}

In the final experiment, we show that we can create a traveling current density perturbation that is amplified in time. We chose $\xi_1=\xi_0$ as before. The Drude weight perturbation frequency is chosen based on an approximate analysis that is described in \cite{tong-thesis}.  The analysis suggests that by choosing $\omega_1 = \sqrt{D_0(\xi_0+\xi_1)/(2\epsilon)} - \omega_0$, we can produce the desired result. This is shown in Figure \ref{FDM-growing}.

\section{Discussion}
In this work, we studied the propagation of plasmons on a graphene sheet whose Drude weight is both time- and space-dependent. Our interest is in the development of an efficient and an accurate numerical method for simulating such phenomena. Directly solving the governing Maxwell's equation with jump conditions on the graphene sheet presents a challenge. First, we need to specify an appropriate absorbing boundary condition around the computational domain. However, it is unclear what is the appropriate condition near the graphene sheet since the plasmons are highly localized to the sheet. Second, we found that due to the dual dependence of the Drude weight, our attempt to build a time-stepping scheme involves inverting a large matrix at each step.

Our approach is to derive a new partial integro-differential equation for the time-integral of the current density on the graphene sheet. This equation is `lives on the graphene sheet and because of the reduced dimension, it can be treated directly. In the present work, the equation is solved by a finite difference method. After numerically establishing convergence of the scheme, we calculated the EM field for a specific type of time- and space-dependent Drude weight. We were able to show that both traveling and standing perturbational fields can be generated.

Finally, it should be repeated here that this work represents a preliminary investigation into partial integro-differential equations of the type that appear in this study. Little is known about such equations though we are convinced that they appear in other applications.

\subsection*{Acknowledgments}
We thank Tony Low at the University of Minnesota for introducing us to this research area. He provided us with research directions and valuable advice. We are also grateful to Paul Martin at the Colorado School of Mines for checking some of our early calculations and for sharing with us some of his insights.

\bibliographystyle{siamplain}
\bibliography{paper2.bbl}

\end{document}